\documentclass[11pt, a4paper]{article}

\usepackage{fullpage}	

\usepackage{amsmath, amssymb, amsthm}
\usepackage{hyperref}

\usepackage[english]{babel}
\usepackage[latin1]{inputenc}
\usepackage{xy}
\usepackage{graphicx}
\usepackage{xcolor}

\usepackage{verbatim}
\usepackage{enumerate}
\usepackage{url}
\usepackage{subfig}

\usepackage{float}

\usepackage{todonotes}

\newcommand{\deta}{\dot\eta}

\newcommand{\E}{\mathbf{E}}
\newcommand{\eps}{\epsilon}
\newcommand{\N}{\mathbb{N}}
\newcommand{\Ob}{\mathcal{D}}
\newcommand{\ptA}{\textsf{A}}
\newcommand{\ptB}{\textsf{B}}
\newcommand{\ptC}{\textsf{C}}
\newcommand{\R}{\mathbb{R}}

\newcommand{\Sb}{\mathbf{S}}
\newcommand{\T}{\mathbb{T}}
\newcommand{\ud}{\mathrm{d}}

\newcommand{\vphi}{\varphi}

\renewcommand{\O}{\mathrm{O}}
\renewcommand{\S}{\mathrm{S}}

\DeclareMathOperator{\Rep}{Re}
\renewcommand{\Re}{\Rep}

\theoremstyle{definition}

\theoremstyle{remark}

\DeclareMathOperator{\Hom}{Hom}
\DeclareMathOperator{\Het}{Het}
\DeclareMathOperator{\atant}{atan2}

\DeclareMathOperator{\tr}{trace}

\newcommand{\rset}[2]{\left\lbrace\, #1\,\left|\;#2\right.\right\rbrace}

\newcommand{\set}[2]{\rset{#1}{#2}}
\newcommand{\sset}[1]{\left\lbrace #1\right\rbrace}
\newcommand{\tsset}[1]{\big\lbrace #1\big\rbrace}
\newcommand{\abs}[1]{\left|{#1}\right|}

\renewcommand{\u}{\mathrm{u}}
\newcommand{\s}{\mathrm{s}}

\newcommand{\shp}{\shortparallel}

\newcommand{\lpar}{\lambda^{\shp}}
\newcommand{\lperp}{\lambda^{\perp}}

\author{Peter Ashwin\\
Department of Mathematics and Statistics\\
University of Exeter, Exeter EX4 4QF, UK
\and
Christian Bick\\
Department of Mathematics\\
Vrije Universiteit Amsterdam,
Amsterdam, NL
}

\title{Global bifurcations organizing weak chimeras in three symmetrically coupled Kuramoto oscillators with inertia}

\begin{document}

\maketitle

\begin{abstract}
Frequency desynchronized attractors cannot appear in identically coupled symmetric phase oscillators because ``overtaking" of phases cannot occur. 
This restriction no longer applies for more general identically coupled oscillators. Hence, it is interesting to understand precisely how frequency synchrony is lost and how invariant sets such as attracting weak chimeras are generated at torus breakup, where the phase description breaks down. 
Maistrenko et al (2016) found numerical evidence of an organizing center for weak chimeras in a system of $N=3$ coupled identical Kuramoto oscillators with inertia. 
This paper identifies this organizing center and shows that it corresponds to a particular type of non-transverse heteroclinic bifurcation that is generic in the context of symmetry. At this codimension two bifurcation there is a splitting of connecting orbits between the in-phase (fully synchronized) state.
This generates a wide variety of associated bifurcations to weak chimeras. 
We further highlight a second organizing center associated with a codimension two symmetry-breaking heteroclinic connection.
\end{abstract}

\maketitle
\allowdisplaybreaks

\section{Introduction}

Coupled oscillators are of great interest for a variety of applications where they provide tractable models of spatio-temporal pattern generation through interaction between systems with attracting periodic behavior. For weak coupling of~$N$ hyperbolic attracting limit cycle oscillators, a standard normal hyperbolicity argument~\cite{hoppensteadt2012weakly} means there is an attracting invariant $N$-torus containing all nearby attractors. 
For a general system of $N=3$ non-identical oscillators and in this limit of weak coupling bifurcations a wide variety of partial frequency synchronized solutions were found in~\cite{Baesens1991}. These involved considering perturbations of the intrinsic frequencies of the oscillators. It is natural to ask whether these frequency-desynchronized solutions can appear spontaneously, rather than through forced symmetry breaking. 

The general case of $N$~symmetrically coupled phase oscillators was considered for example in~\cite{Ashwin1990,Ashwin2008,Ashwin2016}. A result from~\cite{Ashwin1990} is that for such systems with full permutation symmetry~$\Sb_N$, the low dimensionality means that only frequency synchronized solutions can appear, even though a wide variety of phase dynamics are possible. 

The possibility of breaking of frequency synchrony in networks of oscillators that are still identical is discussed in~\cite{ashwin2015weak} where the idea of a {\em weak chimera} as a symmetry breaking of frequencies is introduced. It is also highlighted in~\cite{ashwin2015weak} that because of the presence of 2-in-phase invariant subspaces that are codimension one in the phase space of phase differences, weak chimeras cannot appear in $\Sb_N$-symmetric phase oscillators. But they can appear in a wide range of less symmetric networks that still have a transitive action of symmetries which mean all oscillators can be viewed as identical. For example, they can be found in oscillator networks with modular structures~\cite{ashwin2015weak} or rings with next-nearest neighbor coupling~\cite{thoubaan2018existence}. Similar structures can be found that relate chaotic dynamics~\cite{bick2016chaotic,Bick2015d}.

Such solutions with broken frequency symmetry can also occur for more general $\Sb_N$-symmetric fully symmetric coupled oscillators: 
For example, chimeras can be found in coupled Landau oscillators with full permutation symmetry~\cite{SethiaSen}. This leads to the interesting question of how the breakdown of an invariant $N$-torus associated with weak coupling can lead to the appearance of solutions with broken frequency synchrony. 
The smallest case where one could plausibly see this is for three fully and identically coupled phase oscillators with $\Sb_3$~symmetry. Indeed, this has been observed in 3~coupled Kuramoto-style oscillators with inertia~\cite{Maistrenko2016}, where a wide variety of periodic and chaotic solutions are found in the case that the damping is not too strong. 
Kuramoto oscillators with inertia have been used to model power grid networks~\cite{Filatrella2008} and show a range of dynamics~\cite{Olmi2015}. 
For strong enough damping this reduces to Kuramoto oscillators with identical frequencies where neither chaos nor breaking of frequency synchrony is possible~\cite{Ashwin1990}. 

In this paper, we exhibit some bifurcations associated with heteroclinic networks that organize the emergence of frequency-synchronized solutions in three coupled Kuramoto oscillators with inertia.
Firstly, we clarify the bifurcation happening at the organizing center identified in~\cite{Maistrenko2016} and named in that work as Point~$\ptB$: 
By developing a bespoke Lin's method~\cite{lin1990using,Homburg2010}, we show numerically that this is a symmetric codimension-two homoclinic bifurcation where there is a tangency of stable and unstable manifolds.
Specifically, we (a)~find the space of bounded solutions of the variational equations and (b)~set up a bifurcation function that is zero precisely when there is a tangency of manifolds.
Numerical bifurcation analysis further sheds light on codimension-one bifurcations related to Point~$\ptB$.
Secondly, we identify a further codimension-two Point~$\ptC$ that corresponds to the existence of a heteroclinic connection between synchrony and a nontrivial equilibrium.
This yields an entire network of heteroclinic connections that allow for complicated patterns of oscillator phases passing one another.
These organizing centers are not only interesting examples of global bifurcations in symmetric systems in their own right but also elucidate how heteroclinic networks can lead to the emergence of weak chimeras in minimal networks of phase oscillators with inertia.

The rest of the paper is organized as follows: 
Section~\ref{sec:model} introduces the model of three symmetric Kuramoto oscillators with inertia, and summarizes basic properties of the phase space. 
In Section~\ref{sec:equil}, we discuss the equilibria of the system, their linear stability, and relevant bifurcations that are relevant for the codimension-two Point~$\ptB$.
Section~\ref{sec:bifs} identifies the bifurcation at Point~$\ptB$ as a symmetric homoclinic tangency and confirms this numerically using Lin's method. 
We turn to the global dynamics in parameter space near Point~$\ptB$ in Section~\ref{sec:global} and identify the codimension-two Point~$\ptC$ as an organizing center. 
This point corresponds to a heteroclinic connection between distinct equilibria.
We conclude with some remarks in Section~\ref{sec:discuss}.

\section{Three globally coupled Kuramoto oscillators with inertia}
\label{sec:model}

In this paper, we analyze the dynamics and bifurcations of $N=3$ globally and identically coupled Kuramoto oscillators with inertia as in~\cite{Maistrenko2016}. Specifically, the phase $\theta_k\in\T$ of oscillator~$k$ evolves according to
\begin{align}\label{eq:KuramInertia}
m\ddot\theta_k + \eps\dot\theta_k &= \omega + \frac{\mu}{N}\sum_{j=1}^{N}\sin(\theta_j-\theta_k-\alpha).
\end{align}
for parameters $m>0$, $\omega$, $\mu$, $\epsilon$ and~$\alpha$. 
The system has a continuous phase-shift symmetry $\theta_k\mapsto\theta_k+\phi$ for $\phi\in\T$.
Without loss of generality, we may assume $\omega=0$ by exploiting the phase shift symmetry and set $m=1$ by rescaling~$\eps$ and~$\mu$.
Note that rescaling time $t\mapsto \nu^{-1}t$ with $\nu>0$ corresponds to a parameter transformation $(\eps, \mu)\mapsto(\nu\eps, \nu^2\mu)$.
We henceforth fix $\epsilon=0.1$ as in~\cite{Maistrenko2016}.


Due to the global and identical coupling, the system~\eqref{eq:KuramInertia} is equivariant with respect to the symmetric group of three elements~$\Sb_{3}$ acting by permuting the indices of the oscillators.
The symmetry implies the existence of invariant sets
\begin{align}\label{eq:subsets}
\E_{12} &= \tsset{\theta_1 = \theta_2, \dot\theta_1 = \dot\theta_2}, &
\E_{23} &= \tsset{\theta_2 = \theta_3, \dot\theta_2 = \dot\theta_3}, &
\E_{31} &= \tsset{\theta_3 = \theta_1, \dot\theta_3 = \dot\theta_1},
\end{align}
which correspond to phase configurations where two oscillators are synchronized in phase and frequency. 
Note that~$\Sb_3$ is generated by the three cycle $\rho=(231)$ and the transposition $\kappa=(23)$; 
observe that~$\rho$ permutes the invariant sets~\eqref{eq:subsets} while $\kappa$~permutes~$\E_{12}$ and~$\E_{13}$ and leaves~$\E_{23}$ invariant. 
The invariant sets intersect in the in-phase solution
\begin{equation}
\O = \E_{12}\cap\E_{23}\cap\E_{13} = \tsset{\theta_1 = \theta_2 = \theta_3, \dot\theta_1 = \dot\theta_2=\dot\theta_3}.
\end{equation}

\subsection{Phase-difference coordinates}
\label{sect:PhasDiffCoord}

\begin{figure}
	\centering
\includegraphics[width=12cm]{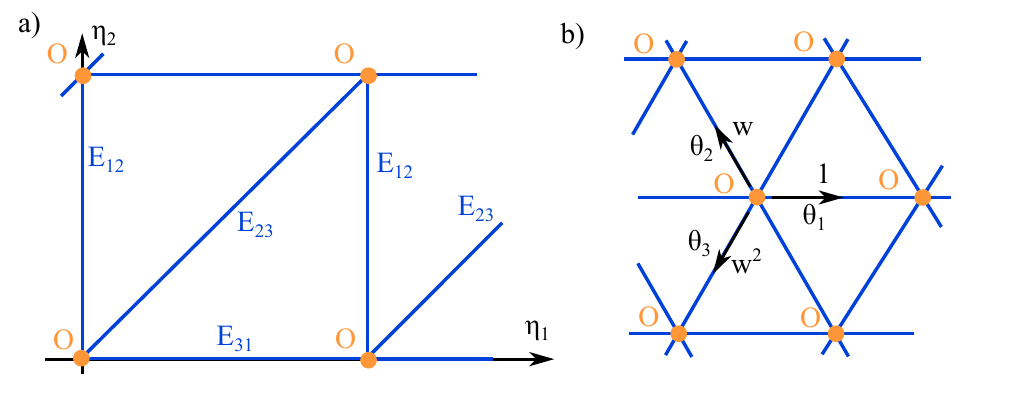}
\caption{\label{fig:CoordinateSystems}
Illustration of coordinates and invariant subspaces for (\ref{eq:KuramInertiaPD}).
Panel~(a) shows a projection of phase space to the phase-difference coordinates~$(\eta_1,\eta_2)$. The invariant subspaces~$\E_{12}$, $\E_{23}$, $\E_{31}$ (blue lines) intersect in the fully synchronized phase configuration~$\O$.
Panel~(b) illustrates the symmetric coordinates~$\varphi$ as defined in~\eqref{eq:SymmCoord}; here the cyclic permutation~$\rho$ of the indices acts as a rotation by~$2\pi/3$.}
\end{figure}

As in~\cite{Maistrenko2016} we introduce phase-difference coordinates $\eta_1 := \theta_1-\theta_2$, $\eta_2:=\theta_1-\theta_3$ to reduce the phase-shift symmetry. 
In these coordinates and setting $m=1$, the system~\eqref{eq:KuramInertia} reduces to
\begin{subequations}
	\label{eq:KuramInertiaPD}
\begin{align}
\ddot\eta_1 + \eps\deta_1 & = -\frac{\mu}{3}\left(2\cos(\alpha)\sin(\eta_1)+\sin(\eta_{2}+\alpha)+\sin(\eta_1-\eta_{2}-\alpha)\right)=: \mu f_1(\eta_1,\eta_2) \\
\ddot\eta_2 + \eps\deta_2 & = -\frac{\mu}{3}\left(2\cos(\alpha)\sin(\eta_2)+\sin(\eta_{1}+\alpha)+\sin(\eta_2-\eta_{1}-\alpha)\right) =: \mu f_2(\eta_1,\eta_2) .
\end{align}
\end{subequations}
that determines the evolution of $\eta := (\eta_1,\dot\eta_1, \eta_2, \dot\eta_2)$.
Writing $\psi_1 := \dot\eta_1$ and $\psi_2:=\dot\eta_2$, \eqref{eq:KuramInertiaPD}~is equivalent to
\begin{subequations}
	\label{eq:KuramInertiaPD4}
	\begin{align}
	\dot{\eta}_1 & = {\psi}_1\\
	\dot\psi_1 &= -\eps\deta_1 -\frac{\mu}{3}\left(2\cos(\alpha)\sin(\eta_1)+\sin(\eta_{2}+\alpha)+\sin(\eta_1-\eta_{2}-\alpha)\right) \\
	\dot{\eta}_2 & = {\psi}_2\\
	\dot\psi_2 &= - \eps\deta_2  -\frac{\mu}{3}\left(2\cos(\alpha)\sin(\eta_2)+\sin(\eta_{1}+\alpha)+\sin(\eta_2-\eta_{1}-\alpha)\right).
	\end{align}
\end{subequations}
In the phase-difference coordinates, the invariant subspaces are given by
\begin{align}
\E_{12} &= \sset{\eta_1 = \dot\eta_1 = 0},&
\E_{23} &= \sset{\eta_1-\eta_2 = \dot\eta_1-\dot\eta_2 = 0},& 
\E_{31} &= \sset{\eta_2 = \dot\eta_2 = 0},
\end{align}
as illustrated in Figure~\ref{fig:CoordinateSystems}(a).
The action of the generators of~$\Sb_{3}$ is given by
$\rho:(\eta_1,\dot\eta_1,\eta_2,\dot\eta_2)\mapsto (-\eta_2,-\dot\eta_1,\eta_1-\eta_2,\dot\eta_1-\dot\eta_2)$ and $\kappa(\eta_1,\psi_1,\eta_2,\psi_2)\mapsto(\eta_2,\psi_2,\eta_1,\psi_1)$.

\subsection{Symmetric coordinates}
\label{sect:SymmetricCoord}

A coordinate system where permutations of the oscillators act orthogonally is given by a complex phase
\begin{equation}
\label{eq:SymmCoord}
\vphi = \theta_1w+\theta_2w^2+\theta_3
\end{equation}
where $w := \exp(2\pi i/3)$; see also~\cite{Ashwin1990}.
These coordinates also reduce the phase-shift symmetry since~$\vphi$ is invariant with respect to it. 
In the complex phase coordinates we have
\begin{align}
    \E_{12} &= \R, & \E_{23} &= \R w, &\E_{31} &= \R w^2
\end{align}
as illustrated in Figure~\ref{fig:CoordinateSystems}(b). Now~$\rho$ acts by multiplication with~$w$. Note that phase-difference coordinates and the complex phase variables are related through
\begin{equation}
(\vphi,\dot\vphi) = (-\eta_1w^2-\eta_2, -\deta_1w^2-\deta_2).
\end{equation}

\section{Equilibria and bifurcations}
\label{sec:equil}

We first give an overview of the equilibria and their stability, and identify the essential bifurcations that organize the dynamics. In Section~\ref{sec:global} we complement the results given in~\cite{Maistrenko2016} by numerical continuation analysis in \textsc{Auto}~\cite{Doedel2000} and by numerical exploration.

Up to symmetry, we classify two equilibria.
First, the in-phase solution~$\O$ at the origin is an equilibrium in reduced coordinates. 
Second, there is one nontrivial equilibrium in each invariant subspace~$\E_{kj}$:
The subspace~$\E_{12}$ contains the equilibrium $\S_{12} = (0, s, 0, 0)$ with 
\[
s=\atant(-6\tan(\alpha), -9+\tan(\alpha)^2),
\]
independent of $\mu$ and $\epsilon$.
Its symmetric copies in $\E_{23}$ and $\E_{31}$ are $\S_{23} = (-s, -s, 0, 0)$ and $\S_{31} = (s, 0, 0, 0)$, respectively.
In the covering space~$\R^4$ of the phase space $(\T\times\R)^2$, we write
\begin{align}
\S_{kj}^{(p,q)} &\in [0+2p\pi, 2\pi+2p\pi)\times\R\times[0+2q\pi, 2\pi+2q\pi)\times\R,\\
\O^{(p,q)} &= (2p\pi, 2q\pi, 0, 0).
\end{align}
denote the equilibria and their symmetric counterparts.
We omit the indices of the specific subspaces or equilibria and just write~$\E$, $\O$, or $\S$ unless the context means the indices are relevant.

\subsection{Linear stability}
\label{sect:Stab}

The structure of the equations of motion and the symmetries restrict the stability properties of the vector field. Let
$\deta =F(\eta,\alpha,\mu)$
denote the vector field of the system of first-order ODEs corresponding to~\eqref{eq:KuramInertiaPD}. 
The linearization of~$F$ at a point $(\eta_1, \psi_1, \eta_2, \psi_2)$ is
\begin{equation}\label{eq:Jacobian}
\ud F = \left(\begin{array}{cccc}
0 & 1 & 0 & 0\\
\mu\partial_{\eta_1} f_1 & -\eps & \mu\partial_{\eta_2}f_1 & 0\\
0 & 0 & 0 & 1\\
\mu \partial_{\eta_1} f_2 & 0 & \mu \partial_{\eta_2}f_2 & -\eps\\
\end{array}\right).
\end{equation}
If~$\lambda_p$ denote the eigenvalues of~$\ud F$, we have 
\begin{equation}
\sum_{p=1}^{4}\lambda_p = \tr(\ud F) = -2\eps.
\end{equation}
For points on the invariant subspaces~$\E_{kj}$ two eigenvalues correspond to eigenvectors that are parallel to~$\E_{kj}$ and two that are transverse; we denote them by $\lambda^\shp_p$ and $\lambda^\perp_p$, respectively. Restricting to the invariant subspaces~$\E_{kj}$ we have that $\sum_{p=1}^{2}\lambda^{\shp}_p = \tr(\ud F|_{\E_{kj}}) = -\eps$ and thus 
\begin{equation}\label{eq:EigRestr}
\lambda^{\shp}_1+\lpar_2=\lperp_1+\lperp_2=-\eps.
\end{equation}

These conditions hold in particular for the linearization at the equilibria~$\O$ and~$\S$ (corresponding to limit cycles in the original system) and restrict their eigenvalues. If $\lambda = a\pm bi$ with $b\neq 0$ is an eigenvalue then~\eqref{eq:EigRestr} implies that $a=-\frac{\eps}{2}<0$. For the equilibrium~$\O$, symmetry implies~$\lpar_p = \lperp_p$.

\subsection{Dynamics and bifurcations of cluster states}
\label{sec:DynamicsInE}

For $\alpha<\frac{\pi}{2}$ the equilibrium~$\O$ is a stable bifocus. 
For $\alpha=\frac{\pi}{2}$ the equilibria~$\O$ and~$\S_{kj}$ undergo a (symmetric) transcritical bifurcation~(TC).
Henceforth we primarily focus on  $\alpha>\frac{\pi}{2}$ beyond the transcritical bifurcation point.
In terms of notation, for equilibria~$P$, $Q$ we write $\Hom(P\to Q)$ for a homoclinic connection from~$P$ to~$Q$; heteroclinic connections are denoted with $\Het$ and $\overset{\E}{\to}$ indicates that the connection lies in~$\E$.

\begin{figure}
\centerline{\includegraphics[width=\linewidth]{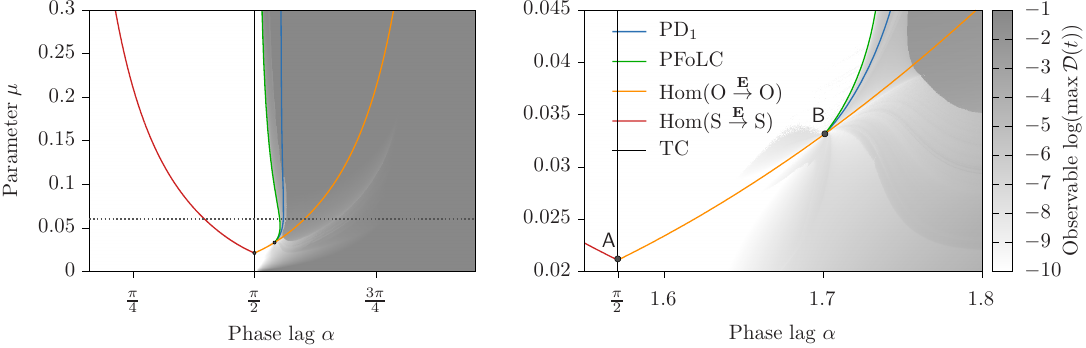}}
\caption{\label{fig:ObsBif}
A bifurcation analysis shows asymptotic dynamics of (\ref{eq:KuramInertia}) with $m=1$, $\epsilon=0.1$ both on and off the invariant set~$\E$; shading represents the observable~$\O$ obtained from numerical integration over~$T=4000$ time units.
Keeping $\mu=0.06$ fixed and varying~$\alpha$ (dashed line in the left panel) there is homoclinic connection $\S\to\S$ within~$\E$ (red line) for $\alpha<\frac{\pi}{2}$ and a homoclinic connection $\O\to\O$ within~$\E$ (orange line) for $\alpha>\frac{\pi}{2}$
that come together at Point~$\ptA$ on the transcritical bifurcation of $\O$~and~$\S$ at $\alpha=\frac{\pi}{2}$ (TC; black line).
A branch of pitchfork of limit cycles (PFoLC) bifurcations on~$\E$ and a secondary period-doubling bifurcation (PD$_1$) come together at the codimension-two Point~$\ptB$ at $(\alpha, \mu)\approx(1.70111, 0.03317)$ that organizes bifurcations to frequency-unlocked solutions.
}
\end{figure}

We first restrict our attention to the dynamics within the invariant (cluster) space~$\E$.
For a fixed~$\mu=0.06$ and varying phase lag~$\alpha$ (indicated by a dashed line in Figure~\ref{fig:ObsBif}).
The saddle~$\O$ and sink~$\S$ have a codimension zero connection $\Het(\O\to \S)$ in~$\E_{kj}$.
A branch of periodic orbits corresponds to partially synchronized weak chimera dynamics: Two oscillators are phase-synchronized and rotate relative to the remaining oscillator.
As~$\alpha$ is decreased, continuation of this branch of periodic orbits leads to a homoclinic bifurcation involving a homoclinic orbit $\Hom\!\big(\S^{(1,0)}\overset{\E}{\to}\S^{(0,0)}\big)$ at $\alpha\approx 1.25$ as~$\alpha$ is decreased (red line in Figure~\ref{fig:ObsBif}).
As~$\alpha$ is increased, the periodic orbit bifurcates in a homoclinic bifurcation involving $\Hom\!\big(\O^{(1,0)}\overset{\E}{\to}\O^{(0,0)}\big)$ at $\alpha\approx 1.9$ (orange line in Figure~\ref{fig:ObsBif}).

Between the red and orange lines shown in Figure~\ref{fig:ObsBif} there is a limit cycle within~$\E$.
The existence of this periodic solution is bounded by the homoclinic bifurcation branches $\Hom\!\big(\O^{(1,0)}\overset{\E}{\to}\O^{(0,0)}\big)$ and $\Hom\!\big(\S^{(1,0)}\overset{\E}{\to}\S^{(0,0)}\big)$.
Continued in two parameters, these branches come together on the transcritical bifurcation line---marked as Point~$\ptA$ in Figure~\ref{fig:ObsBif}.

\subsection{Solutions breaking cluster states}
\label{sec:DynamicsOffE}

Further bifurcations lead to stable dynamics off the invariant set~$\E$ where the oscillators form clusters. 
Continuation in one parameter with~$\mu=0.06$ fixed shows that the branch of periodic orbits in~$\E$ considered above loses stability in a pitchfork bifurcation of limit cycles (PFoLC) at $\alpha\approx1.74$ (green line in Figure~\ref{fig:ObsBif}). 
This results in two symmetry-related branches of a stable limit cycle with partial frequency synchrony off~$\E$ (but broken phase-synchrony\footnote{This transition from a phase- and frequency synchronized weak chimera solution to a solution that synchronized in frequency but not in phase is referred in~\cite{Maistrenko2016} to ``the chimeras becoming imperfect''.}).
These branches lose stability in a period-doubling bifurcation~PD$_1$ (blue line in Figure~\ref{fig:ObsBif}).
Further bifurcations of the periodic orbits involve homoclinic orbits $\Het\!\big(\S_{23}^{(1,0)}\to \S_{23}^{(0,0)}\big)$ and $\Hom\!\big(\S_{23}^{(2,0)}\to \S_{23}^{(0,0)}\big)$ that do not lie in~$\E$ (not shown in Figure~\ref{fig:ObsBif}).

The numerical continuation of these bifurcations in two parameters highlights a codimension-two point that organizes the emergence of these bifurcations between these phase-/frequency-synchrony patterns.
As shown in Figure~\ref{fig:ObsBif}, the PFoLC and PD$_1$ bifurcation curves continued in~$(\alpha,\mu)$ terminate on the point $(\alpha, \mu)\approx(1.70111, 0.03317)$ on the homoclinic bifurcation curve $\Hom\!\big(\O^{(1,0)}\overset{\E}{\to}\O^{(0,0)}\big)$.
This codimension-two point was already observed in~\cite{Maistrenko2016} and dubbed Point~$\ptB$. 
It relates the bifurcations of weak chimera solutions on and off~$\E$.

Indeed, for $\alpha>\frac{\pi}{2}$, the attracting dynamics are typically away from~$\E$.
This is captured by the quantity
\begin{align}
\Ob(\eta) &:= \left(\sin\left(\frac{\eta_1}{2}\right)^2+\deta_1^2\right)\left(\sin\left(\frac{\eta_2}{2}\right)^2+\deta_2^2\right)\left(\sin\left(\frac{\eta_1-\eta_2}{2}\right)^2+(\deta_1-\deta_2)^2\right)
\label{eq:obs}
\end{align}
that vanishes for $\eta\in\E$. 
Thus, the average of~$\Ob(\eta)$ will converge to zero for attractors that lie in~$\E$ while it will be nonzero otherwise.
Note that the dynamics off~$\E$ may be chaotic as the invariant subspace~$\E$ is two-dimensional.
Figure~\ref{fig:ObsBif} shows~$\Ob$ averaged over trajectories for varying parameters together with bifurcation curves computed in \textsc{Auto}, illustrating some of the fine structure of the dynamics near Point~$\ptB$.

\section{Bifurcations of the heteroclinic connection at Point~$\ptB$}
\label{sec:bifs}

We now consider the bifurcation occurring at Point~$\ptB$ in more detail. 
We demonstrate that this is a symmetric homoclinic tangency bifurcation of the $\O\to \O$ homoclinic connection; see~\cite{Homburg2010} for general theory. 
To ease notation we write~$\O$ instead of the specific~$\O^{(p,q)}$ in the covering space unless necessary.

\subsection{A codimension two symmetric homoclinic tangency bifurcation}

First, note that in a large neighbourhood of Point~$\ptB$, the Jacobian~\eqref{eq:Jacobian} of~\eqref{eq:KuramInertiaPD4} near the equilibrium~$\O$ has two identical real positive eigenvalues and two identical real negative eigenvalues.
If $W^\s$~denotes the stable manifold and $W^\u$~the unstable manifold defined in the usual way, we have
\[\dim(W^\u(\O))=\dim(W^\s(\O))=2,\]
i.e., both are two-dimensional invariant manifolds within~$\R^4$. 
To find homoclinic connections $\O\to\O$ we focus on a Poincar\'e section perpendicular to one of the invariant subspaces in~$\E$. For concreteness, we consider the section
\[
\Sigma := \set{(\eta_1,\psi_1,\eta_2,\psi_2)}{\eta_1=-\pi}
\]
perpendicular to~$\E' = \E_{31}$. Note that $\E'$ is the fixed point subspace of $\kappa' =(13)$ that acts on the reduced phase difference coordinates by $\kappa'(\eta_1,\psi_1,\eta_2,\psi_2)\mapsto(\eta_1,\psi_1,-\eta_2,-\psi_2)$

\begin{figure}
    \centering
    \includegraphics[width=0.8\linewidth]
    {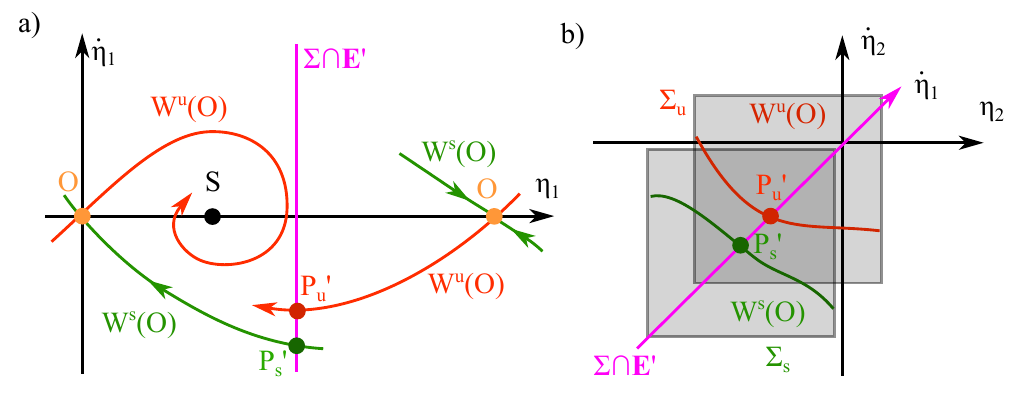}
    \caption{Illustration of the section~$\Sigma$. Panel~(a) shows the dynamics restricted to the invariant subspace~$\E'$: The stable (green) and unstable (red) manifold of~$\O$ intersect $\Sigma\cap\E'$ (purple line) in points~$P_\s'$ (dark green) and $P_\u'$  (red), respectively.
    Panel~(b) shows the section~$\Sigma$:
    The reflection symmetry $\kappa'(\eta_2, \dot\eta_2) = (-\eta_2, -\dot\eta_2)$ that fixes $\Sigma\cap\E'$ forces the intersections of~$W^\u(\O)$ (red line segment) and $W^\s(\O)$ (green line segment) with~$\Sigma$ to affine planes $\Sigma_{\s,\u}$ determined by~$\dot{\eta}_1$ corresponding to $P_{\s,\u}'$ (grey). Observe the codimension zero connection from~$\O$ to~$\S$.
    }
    \label{fig:CodimTwo}
\end{figure}

Restricted to~$\E'$, the section $\Sigma' := \Sigma\cap\E'$ corresponds to a line; cf.~Figure~\ref{fig:CodimTwo}(a). 
Both the unstable manifold~$W^\u(\O)$ and the stable manifold~$W^\s(\O)$ will intersect~$\Sigma'$ in a single point, say~$P'_\u$ and~$P'_\s$, respectively. 
Hence, it is generically codimension one that there are homoclinic connection~$\Gamma$ from $\O\to \O$ within~$\E'$ if $P'=P'_\u=P'_\s = (\pi,p',0,0)$\footnote{It is also generically codimension one that there can be connections that are \emph{not} within~$\E'$.}.

\newcommand{\Tg}{\mathrm{T}}

Now consider the full (reduced four-dimensional) phase space at~$P'\in \Gamma\subset \E'$.
The tangent space at~$P'$ can be decomposed as
\begin{equation}\label{eq:Decomp}
    \R^4 = W_1\oplus W_2 \oplus W_3 \oplus W_4
\end{equation}
where $W_1=\Tg_{P'}W^\s(\O)\cap \Tg_{P'} W^\u(\O)$, $\Tg_{P'}W^\s(\O)=W_1\oplus W_2$ and $\Tg_{P'}W^\u(\O)=W_1\oplus W_3$. 
Given that~$W^\s(\O)$ and~$W^\u(\O)$ are both two-dimensional, it might be expected in a generic situation that 
\begin{align}
\label{eq:dimZ2}
\dim(W_1)&=2
\end{align}
will only appear at codimension three in parameter space:
$W^\s(\O)\cap \Sigma$ and~$W^\u(\O)\cap \Sigma$ close to~$P'$ are line segments and it generically takes a two parameter variation to align them.

Here the problem is restricted by symmetry so that a tangency only requires varying one additional parameter leading to a codimension two homoclinic tangency.
More specifically, note that since~$\E'$---and in particular~$P'$---is fixed under~$\kappa'$ we have that $W^\u(\O)\cap\Sigma$ and $W^\s(\O)\cap\Sigma$ close to~$P'$ into itself. 
This can only happen if $W^\u(\O)\cap\Sigma$, $W^\s(\O)\cap\Sigma\subset \{\dot\eta_1 = p'\}\cap \Sigma$; this is illustrated in Figure~\ref{fig:CodimTwo}(b).
This symmetry-induced constraint reduces the codimension by one allowing for a homoclinic tangency at codimension two.

\subsection{A Lin's method to find point $\ptB$}

\begin{figure}
	\centering
 \includegraphics[width=8cm]{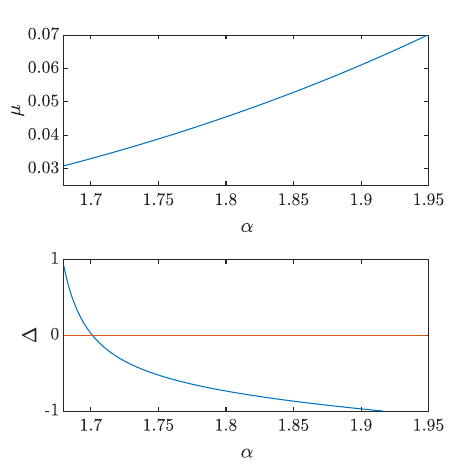}~
 \includegraphics[width=8cm]{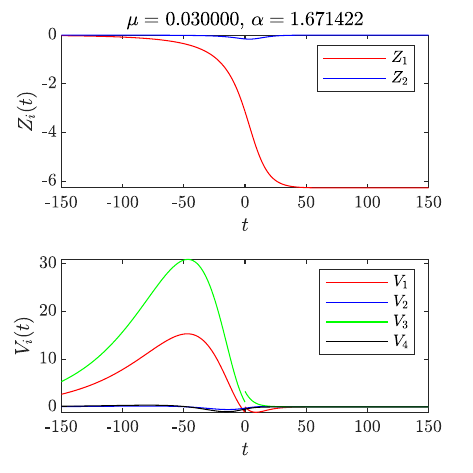}
\caption{The boundary value problems (\ref{eq:BVPhom},\ref{eq:BVPhomvar}) are solved on the truncated interval $[-T,T]$ with $T=150$. Top left: the location $\mu=\mu_{\Gamma}(\alpha)$ of the homoclinic orbit $\Gamma$ to $\O$ as a function of~$\alpha$ near Point~$\ptB$. Bottom left: the ``gap function'' $\Delta(\alpha)$: points where this crosses zero corresponds to a homoclinic $\Gamma$ that satisfies the extra tangency condition of manifolds (\ref{eq:dimZ2}). Note that crossing at $\alpha=1.7011$ and $\mu=0.03317$ at point $\ptB$, indicating the presence of a two-dimensional intersection between the tangent spaces $\Tg_pW^\u(\O)$ and $\Tg_pW^\s(\O)$ for a point on the connection $\Gamma$ homoclinic to~$\O$. Right: the solutions of the boundary value problem for $Z$ (top) and $V$ (bottom) for the initial point away from point $B$: we choose $\mu=0.03$ and use the continuity of $Z$ to determine $\alpha$. Note that $Z$ is continuous meaning we have located the homoclinic but there are discontinuities in the first and third component of $V$ meaning $\Delta\neq 0$.}
\label{fig:lins}
\end{figure}

We now develop a procedure based on Lin's method \cite{lin1990using} to determine the bifurcation point by finding and analysing the connecting orbit; see also~\cite{Homburg2010}.
The first step is to find the homoclinic orbit~$\Gamma$ from~$\O^{(0,0)}$ to~$\O^{(-1,0)}$ in~$\E$, that is,
\[
\Gamma\subset W^\u(\O^{(0,0)})\cap W^\s(\O^{(-1,0)}) \cap \E.
\]
For the Jacobian~$\ud F$
let $w_{\u{}1,2}$ denote the left eigenvectors corresponding to the positive eigenvalues and~$w_{\s{}1,2}$ the left eigenvectors corresponding to the negative eigenvalues.
The problem can be expressed as a three-point boundary value problem on~$[-T,T]$ for large~$T>0$ with four projection boundary conditions and one phase condition.
Specifically, $\Gamma$~is a solution to the boundary value problem
\begin{subequations}
\label{eq:BVPhom}
\begin{align}
\dot{Z} & = F(Z,\alpha,\mu)\\
0 & = \langle w_{s1},(Z(-T)-\O^{(0,0)})\rangle \\
0 & =\langle w_{\s2},(Z(-T)-\O^{(0,0)})\rangle \\
0 & =\langle w_{\u1},(Z(T)-\O^{(-1,0)})\rangle \\
0 & =\langle w_{\u2},(Z(T)-\O^{(-1,0)})\rangle \\
0 & =\langle (Z(0)-(-\pi,0,0,0)),(1,0,0,0)\rangle.
\end{align}
\end{subequations}
We consider $Z=(\eta_1,\dot{\eta}_1,\eta_2,\dot{\eta}_2)$ but work within the invariant subspace $\eta_2=\dot{\eta}_2=0$ and so we only need to solve for the first two components.
The $E^\s$, $E^\u$ denote the stable and unstable subspaces, respectively, the middle four conditions of~\eqref{eq:BVPhom} are equivalent to
\begin{align*}
    Z(-T)&\in E^\u(\O^{(0,0)}), \\ Z(T)&\in E^\s(\O^{(-1,0)})
\end{align*}
while the final equation in~\eqref{eq:BVPhom} is equivalent to the condition
\[P' := Z(0)\in \Sigma.\]
The system~\eqref{eq:BVPhom} is overdetermined, but fixing~$\mu$ and allowing~$\alpha$ to vary will approximate the codimension one condition $P_\u'=P_\s'=P'$ for the homoclinic orbit. This yields the curve~$\mu_{\Gamma}(\alpha)$ of parameter values shown in~\eqref{fig:lins} of parameter values where a homoclinic connection $\O\to\O$ exists in~$\E$; this curve is identical to the orange curve in Figure~\ref{fig:ObsBif} computed using \textsc{Auto}.

The second step is to determine any tangencies of the invariant manifolds.
Note that~$W_1$ in~\eqref{eq:Decomp}---that aligns with the homoclinic connection---can be interpreted as the set of solutions~$V(t)=(\eta_1,\dot{\eta}_1,\eta_2,\dot{\eta}_2)$ of the variational equation $\dot{V}=DF(Z,\alpha,\mu)V$ that remain bounded as $t\rightarrow \pm \infty$.
Thus, in addition to \eqref{eq:BVPhom} we consider solutions of the variational equation for~\ref{eq:BVPhom} given by
\begin{subequations}
	\label{eq:BVPhomvar}
	\begin{align}
	\dot{V} & = DF(Z,\alpha,\mu_h(\alpha))V\\
	0 & = \langle w_{\s1},V(-T)\rangle \\
	0 & =\langle w_{\s2},V(-T)\rangle \\
	0 & =\langle w_{\u1},V(T)\rangle \\
	0 & =\langle w_{\u2},V(T))\rangle \\
	0 & =\langle (1,0,0,0),V(0)\rangle\\
	1 & =\langle (0,0,1,0),V(0)\rangle.
	\end{align}
\end{subequations}
More precisely, we need to find values of $\mu=\mu_h(\alpha)$ that satisfy~\eqref{eq:BVPhom} and where additionally the variational equations have two independent solutions.

We solve this problem using a Lin-type method by allowing solutions within~$V$ to have a discontinuity in two of the components at $t=0$. 
More precisely, we examine solutions of the six-dimensional boundary value problem~\eqref{eq:BVPhom},\eqref{eq:BVPhomvar} in $[-T,T]$ for large~$T$ where all components of~$V$ are continuous at $t=0$ except for the ``gaps''
\begin{subequations}
	\label{eq:BVPhomvapjump}
\begin{align}
\xi_1 &=\langle (1,0,0,0),V(0^+)-V(0^-)\rangle,\\
\xi_2 &=\langle (0,0,1,0),V(0^+)-V(0^-)\rangle,
\end{align}
\end{subequations}
where $0^+$, $0^-$ denote the limits to~$0$ from the positive and negative side, respectively.
Note that since $V(t)=\dot{Z}(t)$ is the only bounded solution of the variational equations that lies within the tangent space to~$\E$ but the final condition on~\eqref{eq:BVPhomvar} ensures that we are finding a nontrivial solution transverse to~$\E$. 
We examine 
\[
\Delta(\alpha):=\xi_2(\alpha,\mu_h(\alpha))
\]
and zeros of this correspond to homoclinic bifurcations that satisfy (\ref{eq:dimZ2}). This is illustrated in Figure~\ref{fig:lins} where we illustrate the $(\alpha,\mu)$ where the homoclinic, the size of the gap $\Delta$ along this curve and a typical solution for the Lin problem. Note that there is a unique location where $\Delta=0$ corresponding to $\ptB$. At this point note that also $\xi_1=0$ because of the symmetry $\kappa$.

\section{Global dynamics, organizing centers, and Point~$\ptC$}
\label{sec:global}

\begin{figure}
    \centering
    \includegraphics[width=15cm]{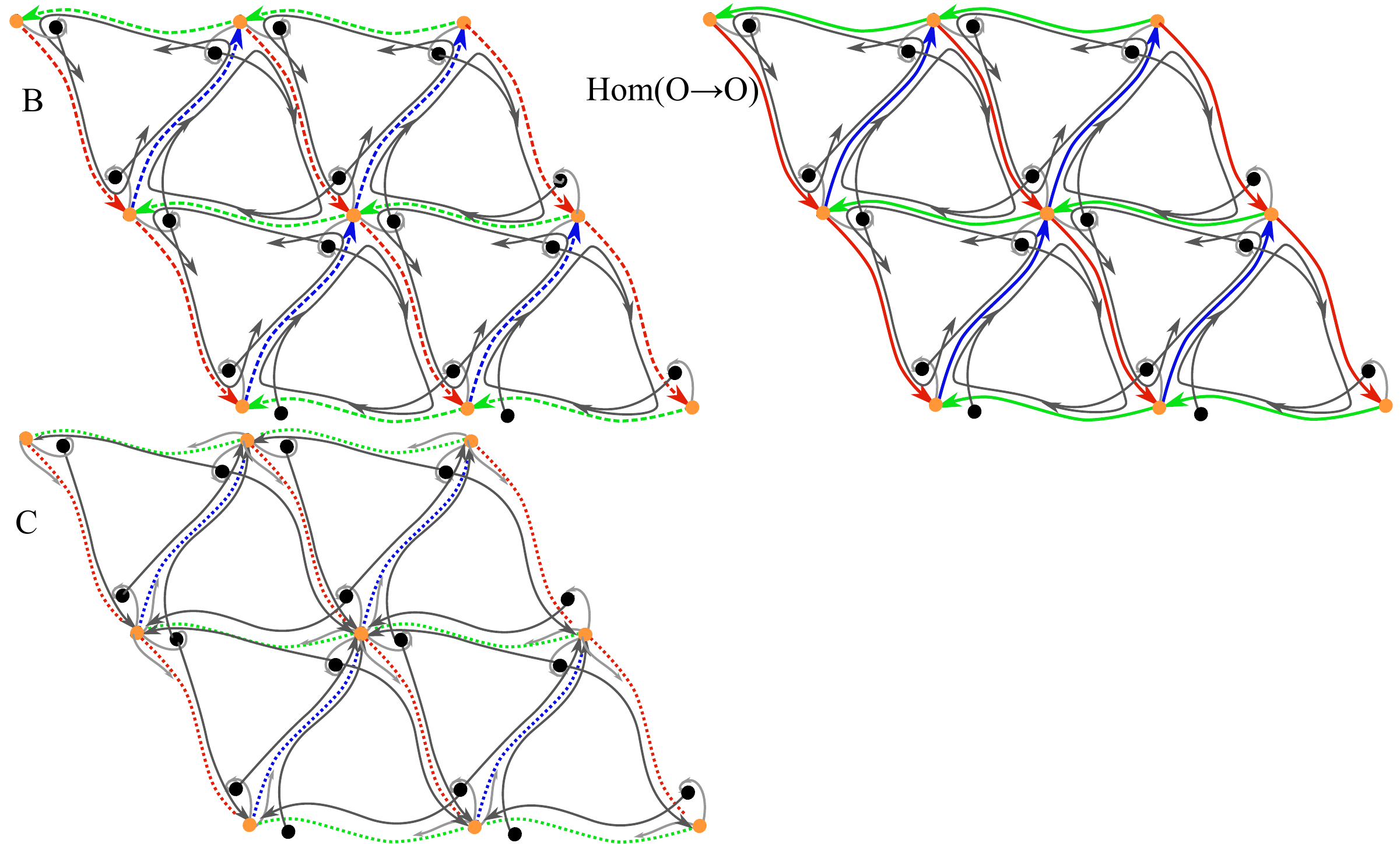}
    \caption{Schematic diagram showing parts of heteroclinic networks connecting the equilibria~$\O$ (orange) and~$\S$ (black) at points~$\ptB$ and~$\ptC$ of the parameter plane as well as for a more general point on the homoclinic bifurcation line $\mbox{Hom}(\O\to \O)$.
    The connections from~$\O$ in $\E_{12}$, $\E_{23}$, $\E_{31}$ are shown in red, green, and blue, respectively. The one-dimensional unstable manifolds of~$\S$ are shown in gray. 
    For~$\ptB$ the $\O\to\O$ connections are shown dashed to indicate they have a degeneracy associated with the two dimensional tangency. The more general $\O\to\O$ homoclinic is as for~$\ptB$ but without the tangency.
    For~$\ptC$ there is a connection from~$\S$ to symmetrically placed copies of~$\O$. There are robust connections from~$\O$ to~$\S$ in one direction and to a saddle periodic orbit (represented by red, green and blue dashed lines). Note that in all of these cases, there will be a variety of pseudo-orbits that wind around the torus while visiting~$\O$ and~$\S$ in various ways.
    }
    \label{fig:colors}
\end{figure}

Varying parameters~$(\alpha,\mu)$ around Point~$\ptB$ shows a wide range of dynamical behavior ranging from periodic to chaotic.
We now show that these dynamics are organized by the interplay of the invariant manifolds of not only~$\O$ but also the nontrivial equilibria~$\S$ that can form heteroclinic connections and networks.
Possibilities are sketched in Figure~\ref{fig:colors}. 
On the one hand, these include not only the codimension-one homoclinic connection $\O\to\O$ in~$\E$ as well as the additional codimension-two tangency at Point~$\ptB$ discussed above.
On the other hand, we find that there is an additional codimension-two Point~$\ptC$ at which there is a heteroclinic $\S\to\O$ connection off~$\E$, which also organizes global dynamics.

\subsection{Heteroclinic networks and chaotic dynamics}

\begin{figure}
		\centering
\centerline{\includegraphics[width=0.5\textwidth]{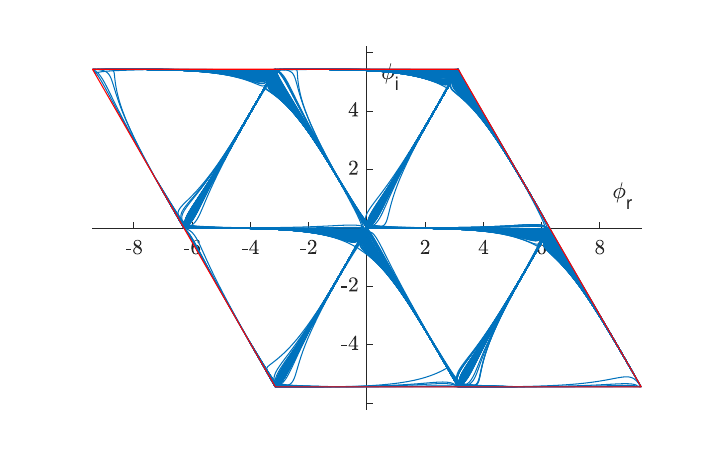}\quad
\includegraphics[width=0.55\textwidth]{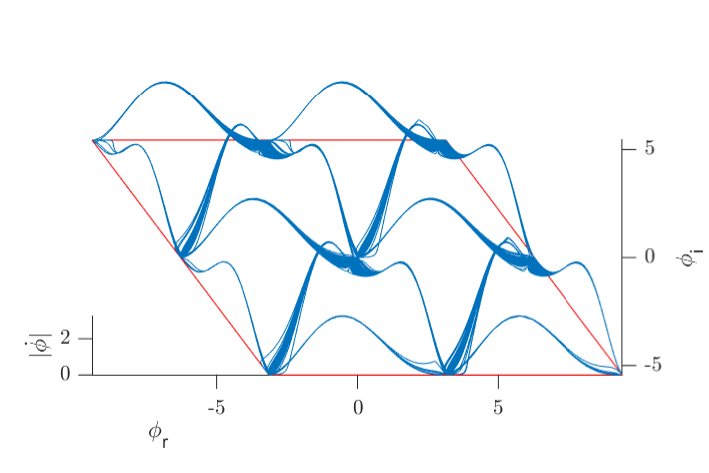}}
\caption{\label{fig:Attractor} Projection of the numerically computed attractor with parameters at Point~$\ptB$ ($\mu=0.03317$, $\alpha = 1.70111$). Left: projection onto $\phi=\phi_r+i\phi_i$ and Right: showing a scaled version of $|\dot{\phi}|^2$ as a third coordinate. The invariant set~$\E$ is depicted by red lines that bound four copies of the unit cell of the $\phi$-torus. The attractor is approximated using a single long trajectory shown in blue - this is given an initial condition close to $\O$ and wraps around the torus in a complicated and recurrent manner.}
\end{figure}

The chaotic attractors that appear close to Point~$\ptB$ are characterized by phase slips associated with chimera states with rapid transitions along the sets~$\E_{kj}$.
Figure~\ref{fig:Attractor} shows a single numerically approximated trajectory started close to~$\O$ that explored the unstable manifold of~$\O$ and hence also the unstable manifolds of~$\S$ in various copies on the torus.
Note that trajectories repeatedly approach~$\S_{kj}$ before being ejected along its one-dimensional unstable manifold. 
The chaotic dynamics near Point~$\ptB$ appear to be organized by the interplay of the invariant manifolds of both~$\O$ and~$\S_{kj}$.
Figure~\ref{fig:colors} illustrates the different heteroclinic connections that are possible:
The codimension-zero $\O\to\S$ connections are shown in purple and the codimension-one $\O\to\O$ connections in~$\E_{12}$, $\E_{23}$, $\E_{31}$ are shown in red, green, and blue, respectively (so colors correspond to phase slips of one of the three oscillators).
Moreover, the one-dimensional unstable manifold of~$\S$ (gray line in Figure~\ref{fig:colors}) approaches~$\O$ and~$\S$ in a complicated sequential manner.

To gain insights into how the unstable manifold~$W^{\mathrm{u}}(\S)$ of~$\S$ wraps around the torus, we follow the sequence of ``phase slips'' along the $\O\to\O$ transition---an increase of phase by an integer multiple of $2\pi$ along one of the invariant subspaces~$\E_{12}, \E_{23},\E_{31}$.
Our approach to visualize the winding of trajectories is inspired by the extraction of kneading sequences from the Lorentz equations in~\cite{Barrio2012}.
More precisely, we assign the $\O\to\O$ connection in each~$\E_{kj}$ a different color (as illustrated in Figure~\ref{fig:colors}).
Colors are then added in the sequence of the phase slips along each of the invariant subspaces~$\E_{kj}$ with a weight that decays according to a geometric sequence:
If~$W^{\mathrm{u}}(\S)$ follows a single~$\E_{kj}$ this results in either red, blue, or green; transitions in the blue and then red direction give purple (with a shift towards red/blue depending on the order).
More details on the algorithm are given in Appendix~\ref{sec:SeqScan}.

\begin{figure}
\centerline{\includegraphics[width=0.8\linewidth]{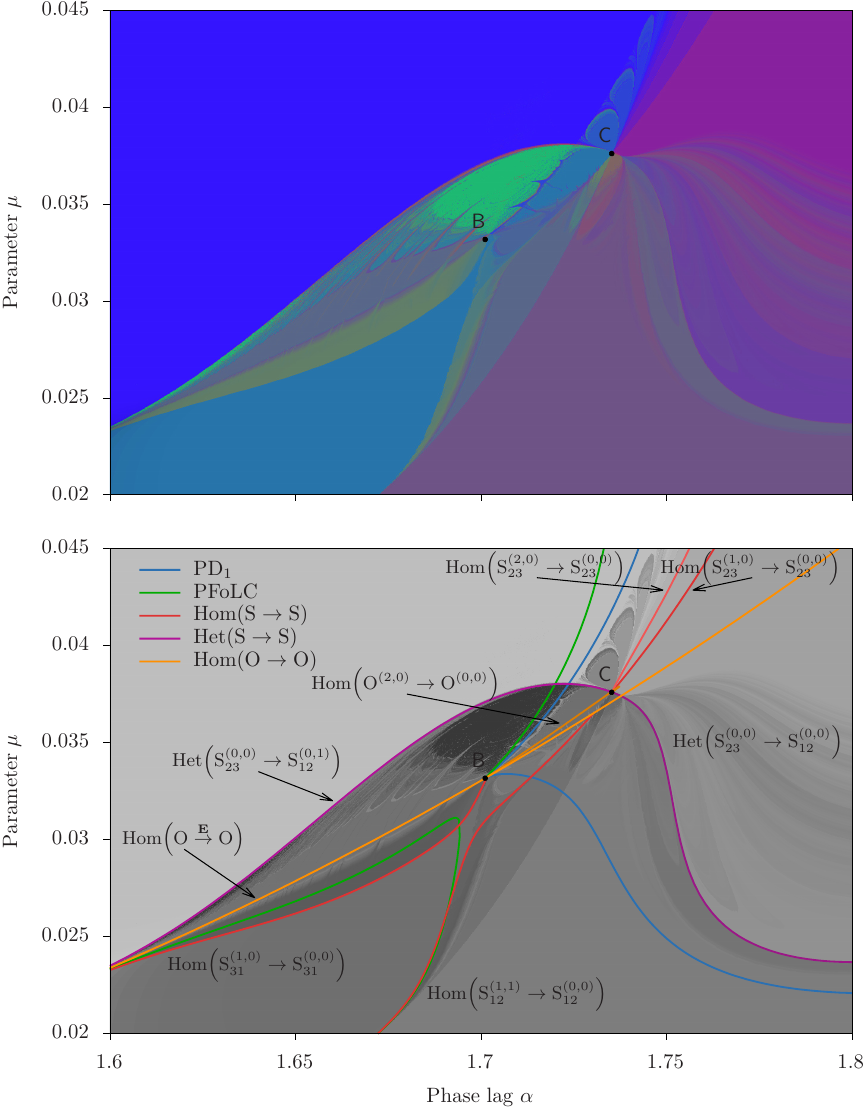}}
\caption{\label{fig:PhaseSlipSeq}Phase-slip sequences obtained by following of~$\S$ capture the global dynamics near Point~$\ptB$  for the system (\ref{eq:KuramInertiaPD4}) with $\epsilon=0.1$.
The left panel shows the colors obtained by integrating the system for $T=4000$ time units as outlined in Appendix~\ref{sec:SeqScan}: Each color corresponds to a phase slip along one of the invariant subspaces~$\E_{kj}$ (cf.~Figure~\ref{fig:colors}); phase slips in different directions give the corresponding mixtures.
The right panel shows two-parameter bifurcation curves obtained through numerical continuation in \textsc{Auto}; colors match those in Figure~\ref{fig:ObsBif} and $\Hom$ indicates a homoclinic connection (modulo~$2\pi$) while $\Het$ indicates a heteroclinic connection between distinct~$\S_{kj}$.
Note that multiple bifurcation curves come together at Points~$\ptB$ with $(\alpha, \mu)\approx(1.70111, 0.03317)$ and~$\ptC$ with $(\alpha,\mu)\approx(1.73521,0.03760)$.
}
\end{figure}

The sequences show the fine structure of the dynamics in parameter space as shown in Figure~\ref{fig:PhaseSlipSeq}.
Distinct colors can be found in the vicinity of Point~$\ptB$ but the picture also unveils a rich structure. 
Further numerical bifurcation analysis (beyond what is shown in Figure~\ref{fig:ObsBif}) shows relevant (local and global) bifurcations that lead to qualitative changes in the phase slip sequences; see Figure~\ref{fig:PhaseSlipSeq}(b).
While Point~$\ptB$ is related to the bifurcation lines shown in Figure~\ref{fig:ObsBif}, there is another codimension-two Point~$\ptC$ that organizes the dynamics and bifurcations.
Inspecting solutions along the bifurcation branches close to Point~$\ptC$ shows that part of the heteroclinic $\O\to\O$ connection off~$\E$ approaches a heteroclinic $\S\to\O$ connection.

\subsection{Heteroclinic connections~$\S\to\O$ at Point~$\ptC$}

Heteroclinic connections $\S\to\O$ are co-dimension two: 
Since the unstable manifold~$W^{\mathrm{u}}(\S)$, it generically intersects a three-dimensional Poincar\'e section~$\Sigma$ in a single point (cf.~Figure~\ref{fig:CodimTwo}).
Variation of two parameters is required for this point to intersect the line~$\Sigma$ that corresponds to the stable manifold of~$\O$.
Note that symmetry forces also heteroclinic connections of the symmetric copies.
Solving the corresponding boundary value problem using XPP confirms the existence of an $\S\to\O$ heteroclinic at 
$(\alpha, \mu) \approx (1.735209263, 0.037605895)$ in agreement of Point~$\ptC$ identified in the previous section.

Note that the $\S\to\O$ connection immediately gives rise to an entire heteroclinic network involving both~$\O$ and all symmetric copies of~$\S$.
Specifically, the codimension-zero connections $\O\to\S_{kj}$ within~$\E$ together with the $\S_{kj}\to\O$ off~$\E$ yields a heteroclinic connection $\O\to\S_{kj}\to\O$.
This means that there are pseudo-orbits that pass by~$\O$ and any symmetric copy of~$\S$ arbitrarily often.
This facilitates the emergence of complicated phase-slip dynamics for parameter values near Point~$\ptC$.

\section{Discussion}
\label{sec:discuss}

We show how certain homoclinic and heteroclinic structures and their bifurcations organize the dynamics in networks of three symmetrically coupled phase oscillators.
First, we clarified the nature of the codimension-two organizing center corresponding to Point~$\ptB$ in the paper of~\cite{Maistrenko2016}:
We find that Point~$\ptB$ is a novel form of symmetric homoclinic tangency bifurcation that can be characterized only using connections between copies of the fully synchronized phase configuration~$\O$. 
The bifurcation will only appear at this relatively low codimension because the presence of the symmetries forces the three invariant subspaces~$\E$ to coincide at~$\O$.
While it was somewhat a surprise that an analysis of connections to and from~$\S$ was \emph{not} necessary to identify the organizing center first described in~\cite{Maistrenko2016}, we identified a second codimension-two Point~$\ptC$ that does involve~$\S$.
The bifurcation at Point~$\ptC$ is easier to identify as it corresponds to a connection between a one-dimensional unstable manifold of~$\S$ and the two-dimensional manifold of~$\O$.

It will be a challenging problem for the future to unfold these bifurcations and to understand the wide variety of behaviors that can be observed to emerge from the organizing centers; cf.~Figure~\ref{fig:PhaseSlipSeq}. 
However, a main challenge to do this can be seen in Figure~\ref{fig:colors}: 
A variety of symmetrically related sections will need to be chosen to ensure that trajectories leaving~$\O$ and~$\S$ are matched appropriately. 
Moreover, the fact that~$\O$ has two-dimensional unstable manifold near~$\ptB$ means that potentially a large number of connections can be involved, including potentially some that connect $\O^{(0,0)}\to\O^{(k,0)}$ for $|k|>1$ as indicated by numerical continuation (cf.~Figure~\ref{fig:PhaseSlipSeq}).
We also expect to find heteroclinic orbits connecting equilibria and periodic orbits that can be continued with an appropriate boundary value setup~\cite{Krauskopf2008a,Doedel2008}.

Together, the codimension-two points identified in this paper organize a wide variety of solutions---many of which partially break frequency synchrony---in its neighborhood. 
As previously noted, these cannot appear if there is an attracting invariant torus and so they are also indicative of crossing the threshold of torus breakdown.
These organizing centers also point towards interesting phenomena of higher codimension: 
For example, at codimension three one would expect the existence of a heteroclinic network consisting of $\S\to\O$ (codimension two), $\O\to\O$ (codimension one), and $\O\to\S$ (codimension zero) connections.

The bifurcations described can presumably be generalized to understand similar organizing centers for systems with larger~$N$, as illustrated in the numerical simulations of~\cite{Maistrenko2016}. 
We have not attempted this but note that the dimensions of the Lin modeling of the variational equation will need to increase correspondingly.

\section*{Acknowledgments}

Thanks to Yuri Maistrenko for suggesting that Point~$\ptB$ was worthy of study. For the purpose of open access, the authors have applied a Creative Commons Attribution (CC--BY) license to any Author Accepted Manuscript version arising from this submission.

\bibliographystyle{unsrt}
\bibliography{refs}

\appendix

\clearpage

\section{Phase-Slip Sequences}
\label{sec:SeqScan}

\newcommand{\vth}{\vartheta}

We extract phase slips of either~$\eta_1$ (along~$\E_{31}$), $\eta_2$ (along~$\E_{12}$), or both~$\eta_1$ and~$\eta_2$ simultaneously (along $\E_{23}$) from the unstable manifold of~$\S$ is one-dimensional close to~$\ptB$.
Note that---compared to evaluating the average frequencies---this strategy captures the transient as well as some indication of the asymptotic dynamics.
For each parameter $(\alpha, \mu)$ we execute the following steps:
\begin{enumerate}
\item Choose an initial condition along the fastest expanding direction of~$\S$: Let~$\lambda_k$ denote the eigenvalues of the linearization of the vector field at~$S$. Suppose that~$\lambda$ is an eigenvalue with corresponding eigenvector~$v$ such that $\Re(\lambda_k)\leq\Re(\lambda)$ for all $k$. For given $\delta>0$, choose $S + \delta\Re(v)$ as initial condition.
\item Solve the system for~$T$ time units.
\item Find times $t_1, t_2, \dotsc$ such that $\xi(t):=\max_{k=1,2}\deta_k^2(t)$ has a local extremum at $t_\ell$ that exceeds $\vth:=\max_{t\in[0, T]}\xi(t)/2$.
\item For each $t_1, t_2, \dotsc$ determine
\begin{subequations}
\begin{align}
\kappa_{12}(m) &:= \begin{cases}1&\text{if } \deta_1(t_m)<0 \text{ and } \deta_1^2(t_m)>\vth, \deta_2^2(t_m)<\vth,\\
-1&\text{if } \deta_1(t_m)>0\text{ and } \deta_1^2(t_m)>\vth, \deta_2^2(t_m)<\vth,\\
0&\text{otherwise}.
\end{cases}\\
\kappa_{23}(m) &:= \begin{cases}1&\text{if } \deta_1(t_m)>0\text{ and } \deta_1^2(t_m)>\vth, \deta_2^2(t_m)>\vth,\\
-1&\text{if } \deta_1(t_m)<0\text{ and } \deta_1^2(t_m)>\vth, \deta_2^2(t_m)>\vth,\\
0&\text{otherwise}.
\end{cases}\\
\kappa_{31}(m) &:= \begin{cases}1&\text{if } \deta_2(t_m)<0\text{ and } \deta_1^2(t_m)<\vth, \deta_2^2(t_m)>\vth,\\
-1&\text{if } \deta_2(t_m)>0\text{ and } \deta_1^2(t_m)<\vth, \deta_2^2(t_m)>\vth,\\
0&\text{otherwise}.
\end{cases}
\end{align}
\end{subequations}
These coefficients correspond to an `orientation' $\O^{(0,0)}\to\O^{(-1,0)}\to\O^{(-1,-1)}\to\O^{(0,0)}$ on~$E$.
\item For given $Q>0$, $M\in\N$ evaluate 
\begin{equation}K_{kj}:=(1-Q)\sum_{m=1}^M \kappa_{kj}(m)Q^m \in [-1, 1].\end{equation}
\end{enumerate}
The sequences encoded by~$K_{kj}$ can be visualized by choosing three linearly independent color vectors $C_\textrm{X}, C_\textrm{Y}, C_\textrm{Z}$ (corresponding to red, green, and blue, for example) by plotting the (appropriately normalized) linear combination
$\abs{K_{12}}C_\textrm{X}+\abs{K_{23}}C_\textrm{Y}+\abs{K_{31}}C_\textrm{Z}$; cf.~Figure~\ref{fig:PhaseSlipSeq}.

\end{document}